# Groupoids Determined by Involutive Automorphisms on Semilattices of Groups

## R. A. R. Monzo[1]

**Abstract.** We explore properties of groupoids $(S, *)$ of the form $x * y = (\alpha x) \circ y$,

for all $\{x, y\} \subseteq S$, where $(S, \circ)$ is a semigroup and a semilattice of groups.

We require $\alpha$ to be an idempotent-fixed automorphism on $(S, *)$ whose square is the identity map. Several characterizations of such groupoids are proved, using analogies of semigroup-theoretic properties and constructions.

## 1. Introduction

Groupoids $S$ of the form $xy = (\alpha x) * (\beta y)$, for all $\{x, y\} \subseteq S$, where $\alpha$ and $\beta$ are endomorphisms of a groupoid $(S, *)$ are called in the literature "linear over $(S, *)$ with parameters $\alpha$ and $\beta$ [cf. eg. 2]. For example, a groupoid is linear over a semigroup $(S, *)$ that is a semilattice of abelian groups with parameters $\alpha$ and $1_S$ if and only if it is a completely inverse AG**-groupoid, where $1_S$ is the identity mapping on $S$, $\alpha$ is an idempotent-fixed automorphism on $(S, *)$ and $\alpha^2 = 1_S$ [3]. The main result of this paper characterizes groupoids that are linear over a semigroup $(S, *)$ that is a semilattice of (not necessarily abelian) groups, with these same parameters. In general we call a groupoid "determined by $\alpha$" if it is linear over a groupoid $(S, *)$ with parameters $\alpha$ and $1_S$, where $\alpha$ is an automorphism on $(S, *)$ and $\alpha^2 = 1_S$. We now state the main result. Any undefined terms will be clarified in Section 2.

**Main Theorem.** *The following statements are equivalent:*

(M1) *$S$ is determined by a mapping $\alpha \in \mathrm{AUT}_e^2(S, *)$ and $(S, *)$ is a semigroup that is a semilattice of groups;*

(M2) *$S$ is a completely inverse groupoid,*
   *there exists $\alpha \in \mathrm{AUT}^2(S)$ such that for all $a,b,c \in S$, $(ab)c = (\alpha a)(bc)$ and*
   *$\mathrm{E}(S)$ is a semilattice __or__ $(ab)^{-1} = (\alpha b^{-1})(\alpha a^{-1})$ for all $a,b \in S$;*

(M3) *$S$ is a strongly regular groupoid,*
   *there exists $\alpha \in \mathrm{AUT}_e^2(S)$ such that for all $a,b,c \in S$, $(ab)c = (\alpha a)(bc)$ and*
   *$\mathrm{E}(S)$ is a semilattice;*

(M4) *$S$ is a completely inverse, right-Bol groupoid,*
   *$\alpha : a \mapsto a(aa^{-1})$ satisfies $\alpha \in \mathrm{AUT}^2(S)$ and*
   *$\mathrm{E}(S)$ is a semilattice __or__ $(ab)^{-1} = (\alpha b^{-1})(\alpha a^{-1})$ for all $a,b \in S$;*

---
[1]10 Albert Mansions, Crouch Hill, London N8 9RE, United Kingdom, e-mail: bobmonzo@talktalk.net

(M5) *S is a disjoint union of groupoids $S(e)$ ($e \in E$), with E a semilattice. Each $S(e)$ ($e \in E$) is determined by a mapping $\alpha_e \in AUT_e^2 G(e)$, with each $G(e)$ ($e \in E$) a group with identity $e$. For each $e,f \in E$ with $f \geq e$ there exists a homomorphism $\delta_{f,e}: S(f) \to S(ef)$ such that*

(1) *for each $e \in E$, $\delta_{e,e}$ is the identity mapping on $S(e)$*

(2) *for $e,f,g \in E$, $g \geq f \geq e$ implies $\delta_{f,e} \, \delta_{g,f} = \delta_{g,e}$*

(3) *for every $b \in S(f)$ and all $e \in E$ we have $\alpha_e (\delta_{f,e} b) = \delta_{f,e} (\alpha_f b)$ and*

(4) *for $a \in S(e)$ and $b \in S(f)$, $ab = (\delta_{e,ef} a) \, (\delta_{f,ef} b)$.*

(M6) *S is a semilattice E(S) of groupoids $S(e)$ ($e \in E$), where each $S(e)$ is determined by a mapping $\alpha_e \in AUT_e^2 G(e)$ and each $G(e)$ is a group with identity $e$. Also, for all $a,b,c \in S$, $(ab)c = (\alpha a)(bc)$, where $\alpha = \bigcup_{e \in E} \alpha_e$.*

Note that the characterizations M2, M3 and M4 use analogies of semigroup-theoretic properties such as inverseness, regularity etc., along with the identity $(ab)c = (\alpha a)(bc)$. The characterization M5 is expressed in terms of a union-of-subgroupoids construction, inspired by the semigroup construction theorem for semilattices of groups. We describe this construction theorem for semilattices of groups in Section 2, along with some well-known properties of semilattices of groups. Finally, the characterization M6 combines the union-of-subgroupoids construction with the identity $(ab)c = (\alpha a)(bc)$.

In Section 3 we prove that groupoids determined by idempotent-fixed involutory automorphisms on a semigroup $(S, *)$ that is a semilattice of groups, have certain properties and we prove a number of interrelationships amongst these properties. All Lemmas in this Section are used later to prove the Main Theorem, in Section 5. Two Lemmas are proved in Section 4 that are also essential to the proof of the Main Theorem. This section is separate to section 3 because its Lemmas deal with semilattices of groupoids, each of which is determined by an idempotent-fixed, involutive automorphism on a group.

## 2. Preliminary definitions and results

All mappings will be written on the left including the mappings $\delta_{f,e}$ of Lemmas 11, 13 and the Main Theorem. The mappings $\delta_{f,e}$ are often written on the right in previous publications [1, p.128].

We denote a groupoid with underlying set $S$ as either $S$ or $(S, *)$, where the mapping $*: S \times S \to S$ is the groupoid product. As a convention we may denote the groupoid product as $x * y$, $x \circ y$ or $xy$. As for parentheses, we use ( ), (( )), [ ], [[ ]], { }, {{ }} and < >. The fact that we are dealing with groupoids, not semigroups, necessitates such notation and hopefully assists the reader.

Then AUT $(S)$ or AUT $(S, *)$ denotes the collection of automorphisms on $S$, $AUT^2(S, *)$ denotes the collection of mappings $\{ \alpha \in AUT(S, *): \alpha^2 = 1_S \}$, where $1_S$ is the identity mapping on $S$, and $AUT_e^2(S, *)$ denotes all idempotent-fixed members of $AUT^2(S, *)$; that is, by definition, $\alpha e = e$, for $e \in E(S) = \{x \in S : x = x^2\}$. A mapping on a set $S$ is called an involutive mapping if $\alpha^2 = 1_S$. So $AUT^2(S, *)$ is the set of involutive automorphisms on $\{S, *\}$. If a mapping $\beta$ is an isomorphism of a groupoid $S$ and a groupoid $T$ then we write $\beta: S \cong T$ or simply $S \cong T$.

If $S$ is a disjoint union of sets $S(e)$, with $e \in \mathrm{E}$, and $\alpha_e : S(e) \to S(e)$ for each $e \in \mathrm{E}$, then $\alpha = \bigcup \alpha_e \, (e \in \mathrm{E}) : S \to S$ denotes the mapping $\alpha \, a = \alpha_e a$, for $a \in S(e)$.

A semilattice E is a groupoid satisfying the identities $xy = yx$, $x = x^2$ and $(xy)z = x(yz)$; that is, E is a commutative, idempotent semigroup.

**Definition 1.** *$S$ is determined by a mapping* $\alpha \in \mathbf{AUT}^2(S, *)$ if $xy = (\alpha \, x) * y$ for all $x, y \in S$.

**Result 1**. *The following statements are equivalent*:

(1.1) $S$ *is determined by a mapping* $\alpha \in \mathrm{AUT}^2(S, *)$;
(1.2) $(S, *)$ *is determined by a mapping* $\alpha \in \mathrm{AUT}^2(S)$.

Proof. $(1.1 \Rightarrow 1.2)$ By definition, for all $x, y \in S$, $xy = (\alpha \, x) * y$. Therefore, $\alpha \, (xy) = \alpha \, [(\alpha \, x) * y]$ and, since $\alpha \in \mathrm{AUT}^2(S, *)$, $\alpha \, (xy) = \alpha \, [(\alpha \, x) * y] = (\alpha \, \alpha \, x) * y = x * \alpha \, y$. But since $xy = (\alpha \, x) * y$ for all $x, y \in S$, $(\alpha \, x)(\alpha \, y) = (\alpha \, \alpha \, x) * \alpha \, y = x * \alpha \, y = \alpha \, (xy)$. Hence, $\alpha \in \mathrm{AUT}^2(S)$. Clearly, $x * y = (\alpha x) \, y$ and so (1.2) is valid.

$(1.2 \Rightarrow 1.1)$ By definition, for all $x, y \in S$, $x * y = (\alpha \, x) \, y$. Therefore, $\alpha \, (x * y) = \alpha \, [(\alpha \, x) \, y]$ and, since $\alpha \in \mathrm{AUT}^2(S)$, $\alpha \, (x * y) = \alpha \, [(\alpha \, x) \, y] = (\alpha \, \alpha \, x)(\alpha \, y) = x(\alpha \, y)$. But since $x * y = (\alpha \, x) \, y$ for all $x, y \in S$, $\alpha \, (x * y) = x(\alpha \, y) = (\alpha \, x) * (\alpha \, y)$. Hence, $\alpha \in \mathrm{AUT}^2(S, *)$. Clearly, $xy = (\alpha x) * y$ and so (1.1) is valid.

**Definition 2.** $S$ is called a ***generalized right-Bol*** groupoid if $[(xy)z] \, w = x \, [(yz)w]$, for all $x, y, z \in S$.

**Result 2**. *If $S$ is determined by a mapping* $\alpha \in \mathrm{AUT}^2(S, *)$ *then*

(2.1) $(xy)z = (\alpha \, x)(yz)$, *for all $x, y, z \in S$ if and only if* $(S, *)$ *is a semigroup*
(2.2) *if* $(S, *)$ *is a semigroup then $S$ is a generalized right Bol groupoid and*
(2.3) $\alpha : (S, *) \cong S$ *if and only if* $x * y = xy$.

Proof. By definition, for all $x, y, z \in S$ we have $(xy)z = [\alpha \, (\alpha \, x * y)] * z$. Then, since $\alpha \in \mathrm{AUT}^2(S, *)$, $(xy)z = [\alpha \, (\alpha \, x * y)] * z = (x * \alpha \, y) * z$. Then $x * (x * y * z) = (\alpha \, x) \, [(\alpha \, \alpha \, y) \, z] = (\alpha \, x) \, (yz)$. So if $\{S, *\}$ is a semigroup then $(xy)z = (x * \alpha \, y) * z = x * (\alpha \, y * z) = (\alpha \, x)(yz)$. Conversely, $(xy)z = (\alpha \, x)(yz)$, for all $x, y, z \in S$ implies $(x * \alpha \, y) * z = x * (\alpha \, y * z)$ for all $x, y, z \in S$. For any $q \in S$, setting $y = \alpha q$ gives $(x * q) * z = x * (q * z)$ and so $S$ is a semigroup. Hence, (2.1) is valid. Then, using (2.1) we have $[(xy)z] \, w = [(\alpha \, x)(yz)] \, w = (\alpha \, \alpha \, x) \, [(yz)w] = x[(yz)w]$ and so $S$ is a generalized right-Bol groupoid. This proves (2.2).

Suppose $\alpha : (S, *) \cong S$. Then, for any $x, y \in S$, $x * y = \alpha \, \alpha \, x * \alpha \, \alpha \, y$. Since $\alpha \in \mathrm{AUT}^2(S, *)$ and $\alpha : S \cong \{S, *\}$, $x * y = \alpha \, \alpha \, x * \alpha \, \alpha \, y = \alpha \, (\alpha \, x * \alpha \, y) = (\alpha \, \alpha \, x) \, (\alpha \, \alpha \, y) = xy$. Conversely, if $x * y = xy$ then, since by Result 1 $\alpha \in \mathrm{AUT}^2 S$, $\alpha \, (x * y) = \alpha \, (xy) = (\alpha \, x)(\alpha \, y)$ and so $\alpha : (S, *) \cong S$ and (2.3) is valid. ∎

**Definition 3.** $S$ is an ***inverse groupoid*** if for every $x \in S$ there exists a unique element $x^{-1} \in S$ such that $(xx^{-1}) \, x = x$ and $(x^{-1}x) \, x^{-1} = x^{-1}$. An inverse groupoid $S$ is called a ***completely inverse groupoid*** if $xx^{-1} = x^{-1}x \in \mathrm{E}(S)$ for all $x \in S$.

**Definition 4.** $S$ is called ***strongly regular*** if for all $a \in S$ there exists $x \in S$ such that $a = (ax)a$ and $ax = xa \in \mathrm{E}(S)$.

**Definition 5.** If E is a semilattice and $\{e, f\} \subseteq \mathrm{E}$ then we define $\boldsymbol{e \leq f}$ if $e = ef \ (= fe)$ and $e < f$ (or $f > e$) if $e \leq f$ and $e \neq f$.

**Definition 6.** We say that the groupoid ***S is a semilattice E = E(S) of sub-groupoids S(e) ($e \in E$)*** if $S$ is a disjoint union of the sub-groupoids $S(e)$ $(e \in \mathrm{E})$ and the product $S(e)S(f) \subseteq S(ef)$ *for all $e, f \in \mathrm{E}$*. In the case that S is a semigroup, each $S(e)$ $(e \in \mathrm{E})$ is a semigroup.

We now list some properties of semigroups that are semilattices of groups, properties well-known to semigroup theorists. These semigroups are also described as unions of groups in which the idempotents commute, or inverse semigroups that are unions of groups [1, page 128]. In these semigroups idempotents are in the centre; that is, $ex = xe$ for every idempotent $e \in \mathrm{E}(S)$ and every $x \in S$ [1, Lemma 4.8]. Since semilattices of groups are inverse semigroups, each element has a unique inverse, with which it commutes. Also, $(xy)^{-1} = y^{-1}x^{-1}$ for every $x, y \in S$. Note also that in a semigroup that is a semilattice E of groups we can identify E with E(S), the semilattice of idempotents of $S$ [1, page 127]. These facts will be used in Section 2.

We now state the structure theorem for semigroups $S$ that are a semilattice E = E(S) of groups for reference later in the paper, as well as for comparison with M5.

**Theorem** [1, Theorem 4.11]. *A semigroup $S$ is a semilattice* E *of groups if and only if there are pairwise disjoint groups $G(e) \big( e \in \mathrm{E} \big)$ and, for each pair of elements $e, f \in \mathrm{E}$ with $e < f$, a homomorphism $\delta_{f, e}: G(f) \to G(ef)$ such that $e < f < g$ implies $\delta_{f,e}\delta_{g,f} = \delta_{g,e}$, $e \in \mathrm{E}$ implies $\delta_{e,e} = 1_{G(e)}$, $S = \bigcup G(e)$ $(e \in E)$ and $a_e b_f = \big( \delta_{e,ef} a_e \big)\big( \delta_{f,ef} b_f \big)$ for all $a_e \in G(e)$ and $b_f \in G(f)$.*

## 3. Groupoids determined by semilattices of groups.

**Lemma 1.** *If $S$ is determined by a mapping $\alpha \in \mathrm{AUT}_e^2(S, *)$ and $(S, *)$ is a semigroup that is a semilattice* E *of groups then:*

(1.1)    *$S$ is a completely inverse groupoid;*

(1.2)    *$\mathrm{E}(S) = \mathrm{E}(S, *)$ is a semilattice and*

(1.3)    *$\alpha \in \mathrm{AUT}_e^2 S$.*

*Furthermore, for all $a, b, c, d \in S$*

(1.4)    *$\alpha\, a = a\,(a^{-1}\,a) = a\,(aa^{-1})$,*

(1.5)    *$(ab)^{-1} = (\alpha\, b^{-1})\,(\alpha\, a^{-1})$,*

(1.6)    *$(ab)c = (\alpha\, a)(bc)$,*

(1.7)    *$a^2\, a^{-1} = a$,*

(1.8)    *$(\alpha\, a)^{-1} = \alpha\, a^{-1}$,*

(1.9)    *$S$ is a generalized right-Bol groupoid; that is, $a\,[(bc)d] = [(ab)c]\,d$,*

(1.10)   *$ea = (\alpha\, a)e$ and $ae = e(\alpha\, a)$, for all $e \in \mathrm{E}(S)$,*

(1.11)   *$e(ab) = (ea)(eb)$, for all $e \in \mathrm{E}(S)$,*

(1.12)   *$(ab)\,(ab)^{-1} = (aa^{-1})\,(bb^{-1}) = (b^{-1}a^{-1})\,(b^{-1}a^{-1})^{-1}$ and*

(1.13)   *$aa^{-1} = (a * a^{-1})$, where $a^{-1}$ is the (unique) inverse of a in $(S, *)$.*

Proof. Let $a \in S$ and define $a^{-1}$ = the (unique) inverse of $a$ in $\{S, *\}$. Then, $[a \, (\alpha \, a^{-1})] \, a = [(\alpha \, a) * (\alpha \, a^{-1})] \, a =$ $= \{\alpha \, [(\alpha \, a) * (\alpha \, a^{-1})]\} * a = (a * a^{-1}) * a = a$. Also, $[(\alpha \, a^{-1}) \, a] \, (\alpha \, a^{-1}) = (a^{-1} * a) * (\alpha \, a^{-1})$ and, since $\alpha \in \mathrm{AUT}_e^2(S, *)$, $[a \, (\alpha \, a^{-1})] \, a = [\alpha \left(a^{-1} * a\right)] * (\alpha \, a^{-1}) = \alpha \, [\left(a^{-1} * a\right) * a^{-1}] = \alpha \, a^{-1}$. Hence, $a$ and $\alpha \, a^{-1}$ are inverses of each other in $S$. We proceed to show that $\alpha \, a^{-1}$ is the unique inverse of $a$ in $S$. That is, $\alpha \, a^{-1} = a^{-1}$.

Suppose that $(ax)a = a$ and $(xa)x = x$. Then $a = (ax)a = \{\alpha \, [(\alpha \, a) * x]\} * a = a * (\alpha \, x) * a$. Hence, $\alpha \, a = (\alpha \, a) * x * (\alpha \, a)$. Also, $x = (xa)x = \{\alpha \, [(\alpha \, x) * a]\} * x = x * (\alpha \, a) * x$. So $x$ and $\alpha \, a$ are inverses of each other in $(S, *)$. Since $\alpha \in \mathrm{AUT}(S, *)$, $\alpha \, x$ and $a$ are inverses of each other in $(S, *)$. That is, $\alpha \, x = a^{-1}$. and so $x = \alpha \, \alpha \, x = \alpha \, a^{-1}$. So $S$ is an inverse groupoid.

Now, using the fact that $a * a^{-1} = a^{-1} * a$ for every $a \in \left(S, *\right)$, we have $a a^{-1} = (\alpha \, a) * (\alpha \, a^{-1}) = \alpha \, (a * a^{-1}) =$ $= a * a^{-1} = a^{-1} * a = (\alpha^2 \, a^{-1}) * a = (\alpha \, a^{-1}) a = a^{-1} \, a$. Finally, $(a a^{-1})(a a^{-1}) = [\alpha \left(a a^{-1}\right)] * (a a^{-1}) =$ $= [\alpha \left(a * a^{-1}\right)] * \left(a * a^{-1}\right)$ and, since $\alpha \in \mathrm{AUT}_e^2(S, *)$ and $a * a^{-1} \in \mathrm{E}\left(S, *\right)$, $(a a^{-1})(a a^{-1}) = \left(a * a^{-1}\right) * \left(a * a^{-1}\right) =$ $= a * a^{-1} = a a^{-1} \in \mathrm{E}\left(S\right)$. Hence, $S$ is a completely inverse groupoid and (1.1) is valid.

Suppose that $e \in \mathrm{E}(S)$. Then $e = e^2 = (\alpha \, e) * e$ and, since $\alpha \in \mathrm{AUT}_e^2 \{S, *\}$, $\alpha \, e = \alpha \, [(\alpha \, e) * e] = e * (\alpha \, e)$. This implies $e = (\alpha \, e) * e = [e * (\alpha \, e)] * e = e * [(\alpha \, e) * e] = e * e \in \mathrm{E}(S, *)$. Hence, $\mathrm{E}(S) \subseteq \mathrm{E}(S, *)$ and, therefore, $\alpha \, e = e$. If $e \in \mathrm{E}(S, *)$ then $\alpha \, e = e$ and $e^2 = (\alpha \, e) * e = e * e = e$ and so $\mathrm{E}(S, *) \subseteq \mathrm{E}(S) \subseteq \mathrm{E}(S, *)$. Therefore, $\mathrm{E}(S) = \mathrm{E}(S, *)$ and $\alpha \, e = e$ for all $e \in \mathrm{E}(S)$. If $\{e, f\} \subseteq \mathrm{E}(S)$ then $ef = (\alpha \, e) * f = e * f = f * e = (\alpha \, f) * e = fe$ and $(ef)g = [\alpha \, (ef)] * g =$ $= [\alpha \, (\alpha \, e * f)] * g = (e * \alpha f) * g = e * (\alpha f * g) = e * (fg) = (\alpha \, e) * (fg) = e(fg)$. Thus, $\mathrm{E}(S)$ is a semilattice and (1.2) is valid. Since we have already proved that $\alpha \, e = e$ for all $e \in \mathrm{E}(S)$, by Result 1, $\alpha \in \mathrm{AUT}_e^2 \, S$ and (1.3) is valid.

Now $a \, (a a^{-1}) = a \, (a^{-1} \, a) = (\alpha \, a) * a^{-1} * a = (\alpha \, a) * \alpha \, (a^{-1} * a) = \alpha \, (a * a^{-1} * a) = \alpha \, a$ and so (1.4) is valid.

Recall that, for $a \in S$, $a^{-1}$ = the (unique) inverse of $a$ in $(S, *)$. Let $x$ be the inverse of $ab$ in $(S, *)$. Then, from the proof of (1.1), $(ab)^{-1} = \alpha \, x$, where $x$ is the inverse of $ab$ in $(S, *)$. Since $ab = (\alpha \, a) * b$, the inverse of $ab$ in $(S, *)$ is equal to $b^{-1} * (\alpha a)^{-1} = b^{-1} * (\alpha a^{-1}) =$ and so the inverse of $ab$ in $S$ is $\alpha \, [b^{-1} * (\alpha a^{-1})] = \left(\alpha b^{-1}\right) * a^{-1}$ $= \left(b^{-1}\right)\left(a^{-1}\right) = \left(\alpha b^{-1}\right)\left(\alpha a^{-1}\right)$, which proves that (1.5) is valid. Note that (1.6) and (1.9) follow from Result 2, (1) and (2) respectively.

Now $a^2 \, a^{-1} = \{\alpha \, [(\alpha \, a) * a]\} * (\alpha \, a^{-1}) = a * (\alpha \, a) * (\alpha \, a^{-1}) = a * [\alpha \, (a * a^{-1})] = a * a * a^{-1} = a * a^{-1} * a = a$, so (1.7) is valid.

Also, $(\alpha \, a)^{-1} = \alpha \, (\alpha \, a^{-1}) = a^{-1} = \alpha \, a^{-1}$ and so (1.8) is valid.

Using the fact that idempotents are central in a semilattice of groups [1, Lemma 4.8], if $e \in \mathrm{E}(S)$ and $a \in S$ then $ea = (\alpha \, e) * a = e * a = a * e = (\alpha \, a)e$. Substituting $\alpha \, a$ for $a$ gives $e(\alpha \, a) = ae$, which proves that (1.10) is valid.

If $e \in E(S)$ and $a,b \in S$ then, using ($6$) and ($10$), $(ea)(eb) = [(\alpha\, a)e]\,[(\alpha\, b)e] = \{[a\,(\alpha\, e)]\,(\alpha\, b)\}\, e =$
$= [(ae)\,(\alpha\, b)]\, e = \{[e\,(\alpha\, a)]\,(\alpha\, b)\}\, e = \{(\alpha\, e)\,[(\alpha\, a)\,(\alpha\, b)]\}\, e = \{e\,[(\alpha\, a)\,(\alpha\, b)]\}\, e =$
$= \{e\,[\alpha\,(ab)]\}\, e = [(ab)\, e] = [\alpha\,(ab)]\, e = e(ab)$, which proves (1.11).

Now using (1.5) and the fact that $(a^{-1})^{-1} = a$ in any inverse groupoid, $(b^{-1}\, a^{-1})^{-1} = (\alpha\, a)(\alpha\, b) = \alpha\,(ab)$.
Now $(ab)(ab)^{-1} = (ab)\,[(\alpha\, b^{-1})\,(\alpha\, a^{-1})] = \alpha\,\{[(ab)\, b^{-1}]\}\,(\alpha\, a^{-1}) = \alpha\,\{[(ab)b^{-1}]\,]a^{-1}\} = \alpha\,\{[\alpha\,(ab)](\,b^{-1}\, a^{-1})\} =$
$= \alpha\,[(\,b^{-1}\, a^{-1})^{-1}\,(\,b^{-1}\, a^{-1})]$. Since, from 1.1, $S$ is a completely inverse groupoid, and from 1.3, $\alpha \in \mathrm{AUT}_e^2\, S$,
$(ab)\,(ab)^{-1} = \alpha\,[(\,b^{-1}\, a^{-1})^{-1}\,(\,b^{-1}\, a^{-1})] = (b^{-1}\, a^{-1})^{-1}\,(b^{-1}\, a^{-1}) = (b^{-1}\, a^{-1})\,(\,b^{-1}\, a^{-1})^{-1}$. Using (1.5) and (1.6),
$(ab)\,(ab)^{-1} = (ab)^{-1}\,(ab) = [\alpha\,(\,b^{-1}\, a^{-1})]\,(ab) = [(b^{-1}\, a^{-1})\, a]\, b = \{[\alpha\,(\,b^{-1})]\,(aa^{-1})\}\, b = b^{-1}\,[(aa^{-1})\, b] =$
$= b^{-1}\,[(aa^{-1})\, b]$, which by (1.10) equals $b^{-1}\,[(\alpha\, b)\,(aa^{-1})] = [(\alpha\, b^{-1})\,(\alpha\, b)]\,(aa^{-1}) = [\alpha\,(b^{-1}b)]\,(aa^{-1}) =$
$= (b^{-1}\, b)\,(aa^{-1}) = (aa^{-1})\,(b^{-1}b) = (aa^{-1})\,(bb^{-1})$, which proves (1.12).

Now we have seen in the proof of (2.1) that $a^{-1} = \alpha\, \mathrm{a}^{-1}$, where $\mathrm{a}^{-1}$ is the (unique) inverse of $a$ in $(S, *)$. Thus,
$aa^{-1} = (\alpha\, a) * (\alpha\, \mathrm{a}^{-1}) = \alpha\,(a * \mathrm{a}^{-1}) = a * \mathrm{a}^{-1}$, which proves (1.13) . ∎

Lemma 1 has established properties of groupoids that are determined by an involutive, idempotent-fixed automorphism on a semigroup that is a semilattice of groups. It will be relied on heavily to prove the Main Theorem in Section 5.

We now prove several inter-relationships amongst the properties (1.1) through (1.13). The following Lemmas will also be used in the proof of the Main Theorem. Lemma 2 establishes conditions under which (1.2) and (1.3) are equivalent.

**Lemma 2.** *If $S$ is an inverse groupoid determined by a mapping $\alpha \in \mathrm{AUT}^2\,(S, *)$, $(S, *)$ is a semigroup and $aa^{-1} \in E(S)$ for all $a \in S$ then $\alpha \in \mathrm{AUT}_e^2\, S$ if and only if $\alpha\, a = a\,(a^{-1}\, a)$ for all $a \in S$.*

Proof. ($\Longrightarrow$) For $a \in S$, $\alpha\, a = [(\alpha\, a)\,(\alpha\, a)^{-1}]\,(\alpha\, a)$. By Result 2, $(xy)z = (\alpha\, x)(yz)$, *for all* $x,y,z \in S$ and so
$\alpha\, a = [(\alpha\, a)\,(\alpha\, a)^{-1}]\,(\alpha\, a) = (\alpha\, a\, \alpha\, a\,)[(\alpha\, a)^{-1}\,(\alpha\, a)] = a\,[\alpha\,(a^{-1}\, a)] = a\,(a^{-1}\, a)$.

($\Longleftarrow$) Let $e \in E(S)$. Then $\alpha\, e = e\,(e^{-1}\, e) = e$. Hence, $\alpha \in \mathrm{AUT}_e^2\, S$. ∎

An inverse groupoid with the generalized right-Bol property and that satisfies (1.3) and (1.5) also satisfies (1.2), as we now prove.

**Lemma 3.** *If $S$ is an inverse, generalized right-Bol groupoid, $\alpha \in \mathrm{AUT}_e^2\, S$ and $(ab)^{-1} = (\alpha\, b^{-1})\,(\alpha\, a^{-1})$ for all $\{a,b\} \subseteq S$ then $E(S)$ is a semilattice.*

Proof. Let $\{e,f,g\} \subseteq E(S)$. Then $(ef)^{-1} = (\alpha\, f^{-1})\,(\alpha\, e^{-1}) = (\alpha\, f)\,(\alpha\, e) = fe$ and so

(1): $$(ef)^{-1} = fe.$$

Since $S$ is a generalized right-Bol groupoid, $[(ef)f]\, f = e[(ff)f] = ef = [(ee)e]\, f = e[(ee)f] = e(ef)$ and so

(2): $$[(ef)f]f = ef = e(ef).$$

Also $(ef)(fe) = (ef)[(ff)e] = \{[(ef)f]f\}e = (ef)e$ and so

(3): $$(ef)(fe) = (ef)e.$$

Also $fe = (ef)^{-1} = [e(ef)]^{-1} = \{\alpha\,[(ef)^{-1}]\}(\alpha\,e^{-1}) = [\alpha\,(fe)](\alpha\,e) = [(\alpha\,f)(\alpha\,e)](\alpha\,e) = (fe)e$ and so

(4): $$fe = (fe)e.$$

Using (1), (2), (3) and (4) we have $(ef)(ef) = (ef)[e(ef)] = (ef)[(ee)(ef)] = \{[(ef)e]e\}(ef) = \{[(ef)(fe)]e\}(ef) = = (ef)\{[(fe)e](ef)\} = (ef)[(fe)(ef)] = \{[(ef)f]e\}(ef) = [(ef)e](ef) = [(ef)(fe)](ef) = ef \in \mathrm{E}(S)$. Hence, $ef = (ef)^{-1} = fe$. Since $[(ef)g] = [g(ef)] = [g(fe)] = g[(ff)e] = [(gf)f]e = (gf)e = e(gf) = e(fg)$, $\mathrm{E}(S)$ is a semilattice. ∎

If an inverse groupoid $S$ satisfies $aa^{-1} \in \mathrm{E}(S)$ for all $a \in S$, and if it has an involutive automorphism $\alpha$ satisfying (1.6), then $S$ satisfies (1.3) and (1.4). This is proved in Lemma 4, which also establishes conditions under which (1.2) is equivalent to (1.5).

**Lemma 4.** *If $S$ is an inverse groupoid, $aa^{-1} \in \mathrm{E}(S)$ for all $a \in S$, $\alpha \in \mathrm{AUT}^2(S)$ and $(ab)c = (\alpha\,a)(bc)$ for all $\{a,b,c\} \in S$ then*

(4.1) $\alpha \in \mathrm{AUT}_e^2(S)$,

(4.2) $\alpha\,a = a\,(a^{-1}\,a)$ *for all $a \in S$ and*

(4.3) $\mathrm{E}(S)$ *is a semilattice if and only if $(ab)^{-1} = (\alpha\,b^{-1})\,(\alpha\,a^{-1})$ for all $a,b \in S$.*

Proof. Let $e \in \mathrm{E}(S)$. Then since $\alpha\,e = [(\alpha\,e)(\alpha\,e)^{-1}](\alpha\,e) = e[(\alpha\,e)^{-1}(\alpha\,e)] = e[(\alpha\,e^{-1})(\alpha\,e)] = e\,\alpha\,(e^{-1}\,e) = = e(\alpha\,e)$. Since $\alpha \in \mathrm{AUT}^2(S)$, $e = (\alpha\,e)e$. Hence, $[e(\alpha\,e)]e = (\alpha\,e)e = e$ and $[(\alpha\,e)e](\alpha\,e) = e(\alpha\,e) = \alpha\,e$. Therefore, $\alpha\,e = e^{-1} = e$ and $\alpha \in \mathrm{AUT}_e^2(S)$, proving (4.1).

Since $\alpha \in \mathrm{AUT}_e^2(S)$, $a\,(a^{-1}\,a) = a[\alpha\,(a^{-1}\,a)] = a[\alpha\,(a^{-1})(\alpha\,a)] = [(\alpha\,a)\,\alpha\,(a^{-1})](\alpha\,a) = \alpha\,[(aa^{-1})a] = = \alpha\,a$ and so (4.2) is valid.

Suppose that $\mathrm{E}(S)$ is a semilattice. Let $a,b \in S$. So, using (4.1) and the fact that $(ab)\,c = (\alpha\,a)\,(bc)$ for all $\{a,b,c\} \in S$,

$\{(ab)\,[\alpha\,(b^{-1}\,a^{-1})]\}\,(ab) = [\alpha\,(ab)]\,\{[\alpha\,(b^{-1}\,a^{-1})]\,(ab)\} = [\alpha\,(ab)]\,\{[(b^{-1}\,a^{-1})\,a]\,b\} = = [\alpha\,(ab)]\,\{[\,(\alpha\,b^{-1})\,(a^{-1}\,a)]\,b\} = [(\alpha\,a)(\alpha\,b)]\,\{b^{-1}\,[(\,a^{-1}\,a)\,b]\} = = a\,[(\alpha\,b)\,\{b^{-1}\,[(\,a^{-1}\,a)\,b]\}] = a\,\{(bb^{-1})\,[(\,a^{-1}\,a)\,b]\} = a\,[\,\{[\alpha\,(bb^{-1})]\,(\,a^{-1}\,a)\}\,b\,] = a\,[[(bb^{-1})\,(a^{-1}\,a)]\,b\} = = a\,\{[(\,a^{-1}\,a)(\,bb^{-1})]\,b\} = a\,\{[\alpha\,(\,a^{-1}\,a)][(\,bb^{-1})\,b]\} = a\{(a^{-1}\,a)\,b] = [(\alpha\,a)\,(\,a^{-1}\,a)]\,b = [(aa^{-1})\,a]\,b = ab$. So

(1): $$\{(ab)\,[\alpha\,(b^{-1}\,a^{-1})]\}\,(ab) = ab$$

Using (4.1), (4.2) and the hypotheses of Lemma 4, we prove that $\{[\alpha\,(b^{-1}\,a^{-1})]\,(ab)\}\,[\alpha\,(b^{-1}\,a^{-1})] = = \alpha\,(b^{-1}\,a^{-1})$. Firstly, $\{[\,\alpha\,(b^{-1}\,a^{-1})]\,(ab)\}\,[\,\alpha\,(b^{-1}\,a^{-1})]=(b^{-1}\,a^{-1})\,\{(ab)\,[\,\alpha\,(b^{-1}\,a^{-1})]\}=$

$=(b^{-1}\,a^{-1})\,\{(ab)\,[(\alpha\,b^{-1})\,(\alpha\,a^{-1})]\}=(b^{-1}\,a^{-1})\,\{\,([\,\alpha\,(ab)]\,[\,\alpha\,b^{-1}])\,(\alpha\,a^{-1})\}=$

$=\{\ \alpha\,(b^{-1}\,a^{-1})\,[\,\{\,[\,\alpha\,(ab)][\,\alpha\,b^{-1}]\,\}\,]\ \}\,(\alpha\,a^{-1})\ =\{\ \alpha\ ([b^{-1}\,a^{-1}]\,[\,(ab)b^{-1}\,])\ \}\,(\alpha\,a^{-1})=$

$=\alpha\ \{\ ([b^{-1}\,a^{-1}]\,[\,(ab)b^{-1}\,])\,a^{-1}\,\}$ , so

(2):  $\qquad\qquad\qquad \{[\,\alpha\,(b^{-1}\,a^{-1})]\,(ab)\}\,[\,\alpha\,(b^{-1}\,a^{-1})]=\alpha\ \{\ ([b^{-1}\,a^{-1}]\,[\,(ab)b^{-1}\,])\,a^{-1}\,\}$

But $([b^{-1}\,a^{-1}]\,[\,(ab)b^{-1}\,])\,a^{-1}=[\,\alpha\,(b^{-1}\,a^{-1})]\,\{[(ab)b^{-1}]a^{-1}\}=[\,\alpha\,(b^{-1}\,a^{-1})]\,\{[\,\alpha\,(ab)](b^{-1}\,a^{-1})\}=$

$=\{(b^{-1}\,a^{-1})\,[\,\alpha\,(ab)]\}\,(b^{-1}\,a^{-1})$ , and so by (1),

(3):  $\qquad\qquad\qquad ([b^{-1}\,a^{-1}]\,[\,(ab)b^{-1}\,])\,a^{-1}=b^{-1}\,a^{-1}.$

So (2) and (3) imply that $\{[\,\alpha\,(b^{-1}\,a^{-1})]\,(ab)\}\,[\,\alpha\,(b^{-1}\,a^{-1})]=\alpha\,(\,b^{-1}\,a^{-1})$, which together with (1) implies that $(ab)^{-1}=(\alpha\,b^{-1})(\alpha\,a^{-1})$ for all $a,b\in S$.

Conversely, if $(ab)^{-1}=(\alpha\,b^{-1})(\alpha\,a^{-1})$ for all $a,b\in S$ then, since $(ab)c=(\alpha\,a)(bc)$ for all $\{a,b,c\}\in S$ implies that $S$ is a generalized right-Bol groupoid and since $\alpha\in\mathrm{AUT}_e^{\,2}(S)$, by Lemma 3, $\mathrm{E}(S)$ is a semilattice. Therefore (4.3) is valid. ∎

The following Lemma gives a characterization of groupoids that are determined by involutive, idempotent-fixed automorphisms on semigroups that are semilattices of groups.

**Lemma 5.** *$S$ is determined by a mapping* $\alpha\in\mathrm{AUT}_e^{\,2}(S,*)$ *and* $(S,*)$ *is semigroup that is a semilattice* E *of groups if and only if*

(5.1) *$S$ is a completely inverse groupoid*,

(5.2) *there exists $\alpha\in\mathrm{AUT}^{\,2}(S)$ such that for all $a,b,c\in S$,* $(ab)c=(\alpha\,a)(bc)$ *and*

(5.3) *$\mathrm{E}(S)$ is a semilattice **<u>or</u>** $(ab)^{-1}=(\alpha\,b^{-1})\,(\alpha\,a^{-1})$ for all $a,b\in S$.*

Proof. ($\Longrightarrow$) (5.1), (5.2) and (5.3) follow from Lemma 1.

($\Longleftarrow$) Assume that (5.1), (5.2) and (5.3) are valid. By Lemma 4, $\mathrm{E}(S)$ is a semilattice **<u>and</u>** $(ab)^{-1}=(\alpha\,b^{-1})\,(\alpha\,a^{-1})$ *for all $a,b\in S$.* We define $(S,*)$ as follows: for $a,b\in S$, $a*b=(\alpha\,a)b$. Then $(a*b)*c=[(\alpha\,a)\,b]*c=[a\,(\alpha\,b)]c=(\alpha\,a)\,[(\alpha\,b)\,c]=a*(b*c)$ and so $(S,*)$ is a semigroup. Note that by Lemma 4, $\alpha\,a=a(a^{-1})$ for all $a\in S$ and so $(a*a)*(\alpha\,a^{-1})=[(\alpha\,a)\,a]*(\alpha\,a^{-1})=[a\,(\alpha\,a)]\,(\alpha\,a^{-1})=$

$=(\alpha\,a)\,[(\alpha\,a)\,(\alpha\,a^{-1})]=(\alpha\,a)\,[(\alpha\,a^{-1})\,(\alpha\,a)]=\alpha\,[a\,(a^{-1}a)]=\alpha\,a=a.$

Also, Lemma 4 implies $\alpha\in\mathrm{AUT}_e^{\,2}(S)$ and so $(\alpha\,a^{-1})*(a*a)=(\alpha\,a^{-1})*[(\alpha\,a)\,a]=a^{-1}\,[(\alpha\,a)\,a]=$

$=[(\alpha\,a^{-1})\,(\alpha\,a)]\,a=[\,\alpha\,(a^{-1}a)]\,a=(a^{-1}\,a)\,a=(aa^{-1})\,a=a.$

So we have proved that $a\in(S*a^2)\bigcap(a^2*S)$. By [1, Theorem 4.3], $(S,*)$ is a union of groups.

Suppose that $e, f \in E(S, *)$. Then $e = e * e = (\alpha e) e$ and so $e^{-1} = [(\alpha e) e]^{-1} = (\alpha e^{-1}) [\alpha (\alpha e)^{-1}] =$
$= (\alpha e^{-1}) [\alpha (\alpha e^{-1})] = (\alpha e^{-1}) e^{-1}$. Then $e = (ee^{-1}) e = \{e [(\alpha e^{-1}) e^{-1}]\} e = \{[(\alpha e) (\alpha e^{-1})] e^{-1}\} e =$
$= \{[\alpha (ee^{-1})] e^{-1}\} e = [(ee^{-1}) e^{-1}] e = [(e^{-1} e) e^{-1}] e = e^{-1} e$ and $e^{-1} = (e^{-1} e) e^{-1} = ee^{-1} = e^{-1} e = e$. Therefore
$ee^{-1} = e^{-1} e = e \in E(S)$. Hence, $\alpha e = e$. So then $e * f = (\alpha e) f = ef = fe = (\alpha f) e = f * e$ and so $E(S, *)$ is a semilattice.
By [1, Theorem 4.11], $(S, *)$ is a semilattice of groups. Also $ab = (\alpha a) * b$. Using Lemma 4, $S$ is determined by
$\alpha \in \text{AUT}_e^2 (S, *)$ and the semigroup $(S, *)$ is a semilattice of groups . ∎

The next Lemma shows that a strongly regular groupoid satisfying (1.2), (1.3) and (1.6) is completely inverse.
This allow us to prove a different characterization in Lemma 7, of groupoids determined by idempotent-fixed,
involutive automorphisms on semigroups that are semilattices of groups, with condition (5.1) weakened from
completely inverse to strongly regular, as long as $\alpha$ is idempotent-fixed.

**Lemma 6.** *Suppose that in a strongly regular groupoid $S$ in which* $E(S)$ *is a semilattice and there exists*
$\alpha \in \text{AUT}_e^2 (S)$ *such that for all $x, y, z \in S$, $(xy)z = (\alpha x)(yz)$. Then $S$ is a completely inverse groupoid.*

Proof. Note that, as in the proof of Result 2, $(xy) z = (\alpha x) (yz)$ for all $x, y, z \in S$ implies that $S$ is a right-Bol
groupoid. For each $a \in S$ there exists $x \in S$ such that $a = (ax) a$ and $ax = xa \in E(S)$. Now, $a [(xa) x] = [(ax) a] x =$
$= ax = xa = x [(ax) a] = [(xa) x] a$. Then $\{a [(xa) x]\} a = (ax) a = a$ and $\{[(xa) x] a\} [(xa) x] = (xa) [(xa) x] =$
$= \{[(xa) x] a\} x = (xa) x$.

So we have proved that $a$ and $(xa) x$ are inverses and that $a [(xa) x] = [(xa) x] a = ax = xa \in E(S)$. We need only
therefore show that $(xa) x$ is the unique inverse of $a$. Suppose that $y$ is an inverse of $a$. We prove that $y = (xa) x$.

We have $(ya) y = y$ and $(ay) a = a$. Now, $ax = [(ay) a] x = [\alpha (ay)] (ax)$ and $ay = [(ax)a]y = [\alpha (ax)](ay) = (ax)(ay)$.
But then $(ay)^2 = (ay) [(ax) (ay)] = \{[(ay) a] x\} (ay) = (ax) (ay) = ay \in E(S)$ and so $ax = [\alpha (ay)] (ax) = (ay) (ax)$.
Since $E(S)$ is a semilattice, $xa = ax = (ay) (ax) = (ax) (ay) = ay$.

Since $S$ is strongly regular we can choose $w \in S$ such that $y = (yw) y$ and $yw = wy \in E(S)$. Then $ya = [(yw) y] a =$
$= [\alpha (yw)] (ya) = (yw) (ya)$ and so $(ya)^2 = (ya) [(yw) (ya)] = \{[\alpha (ya)] (yw)\} (ya) = \{[(ya) y] w\} (ya) = (yw) (ya) =$
$= ya$ and so $ya \in E(S)$ and $\alpha (ya) = ya$. Then $xa = x [(ay) a] = [(xa) y] a = [\alpha (xa)] (ya) = (xa) (ya)$ and
$ya = y [(ax) a] = y [(\alpha a) (xa)] = [(\alpha y) (\alpha a)] (xa) = [\alpha (ya)] (xa) = (ya) (xa) = (xa) (ya) = xa = ax = ay = ya$.
Then, $y = (ya) y = (\alpha y) (ay) = (\alpha y) (ax) = (ya) x = (xa) x$, which is what we needed to prove. ∎

**Lemma 7.** *$S$ is determined by a mapping $\alpha \in \text{AUT}_e^2 (S, *)$ and $(S, *)$ is a semigroup that is a semilattice $E$ of*
*groups if and only if*

(7.1) *$S$ is a strongly regular groupoid*,
(7.2) *there exists $\alpha \in \text{AUT}_e^2 (S)$ such that $(ab)c = (\alpha a)(bc)$ for all $a, b, c \in S$ and*
(7.3) $E(S)$ *is a semilattice*.

Proof. ($\Longrightarrow$) (7.1), (7.2) and (7.3) follow from Theorem 2.

($\Longleftarrow$) Assume that (7.1), (7.2) and (7.3) are valid. By Theorem 6, $S$ is completely inverse. Then, by
Theorem 5, $S$ is determined by a mapping $\alpha \in \text{AUT}_e^2 (S, *)$ and $(S, *)$ is semigroup that is a semilattice E of
groups. ∎

The next Lemma establishes conditions under which (1.6) is equivalent to the generalized right-Bol property.

**Lemma 8.** *If $S$ is a completely inverse groupoid satisfying*

(A) $\alpha : a \mapsto a\,(aa^{-1})$ *satisfies* $\alpha \in \mathrm{AUT}^2(S)$ *and*

(B) $\mathrm{E}(S)$ *is a semilattice* **or** $(ab)^{-1} = (\alpha\,b^{-1})\,(\alpha\,a^{-1})$ *for all $a,b \in S$*

*then $(ab)c = (\alpha\,a)(bc)$ for all $a,b,c \in S$, if and only if $S$ is a generalized right-Bol groupoid.*

Proof. ($\Longrightarrow$) If $(ab)c = (\alpha\,a)(bc)$ for all $a,b,c \in S$ then, as in the proof of Result 2, $S$ is a generalized right-Bol groupoid.

($\Longleftarrow$) Note that for $e \in \mathrm{E}(S)$, $\alpha\,e = e\,(ee^{-1}) = e$ and so $\alpha \in \mathrm{AUT}_e^2(S)$. Also, if $(ab)^{-1} = (\alpha\,b^{-1})(\alpha\,a^{-1})$ for all $a,b \in S$ then, by Lemma 3, $\mathrm{E}(S)$ is a semilattice. So (B) implies that $\mathrm{E}(S)$ is a semilattice.

By hypothesis, $S$ is a generalized right-Bol groupoid and we will use that fact, and the fact that $\mathrm{E}(S)$ is a semilattice, without mention throughout the remainder of the proof. We now prove that $(ab)^{-1} = (\alpha\,b^{-1})\,(\alpha\,a^{-1})$ for all $a,b \in S$. Firstly,

$\{(ab)\,[(\alpha\,b^{-1})(\alpha\,a^{-1})]\}\,(ab)\ =$

$= [\,(ab)\ \{\,[b^{-1}(b^{-1}b)]\,[a^{-1}(a^{-1}a)]\,\}\ ]\ (ab) =$

$= (\ \{\,[(ab)b^{-1}]\,(b^{-1}b)\,\}\ \ [a^{-1}(a^{-1}a)]\ )\ (ab)\ =$

$= [(ab)b^{-1}]\ \ \ [\,\{\,(b^{-1}b)\,[a^{-1}(a^{-1}a)]\,\}\ (ab)\,]\ =$

$= [(ab)b^{-1}]\ \ \ \{\,[\,(b^{-1}b)\ \{\,[(a^{-1}\,a)\,a^{-1}]\,(a^{-1}a)\,\}\,]\ (ab)\,\} =$

$= [(ab)b^{-1}]\ \ \ <\ [[\ \{\,[(b^{-1}b)\,(a^{-1}\,a)]\,a^{-1}\,\}\ \ (a^{-1}a)\,]]\ \ \ (ab)\ >\ =$

$= [(ab)b^{-1}]\ \ \ <\ [\,(b^{-1}b)\,(a^{-1}\,a)]\ \ \{\,[a^{-1}(a^{-1}a)]\ (ab)\,\}\ >\ \mathrm{so}$

(1):   $\{(ab)\,[(\alpha\,b^{-1})\,(\alpha\,a^{-1})]\}\,(ab) = [(ab)b^{-1}]\ \ \ <\ [\,(b^{-1}b)\,(a^{-1}\,a)]\ \ \ \{\,[a^{-1}(a^{-1}a)]\ (ab)\,\}\ >$

Note that $a = \alpha\ \alpha\,a = \alpha\,[a\,(aa^{-1})] = (\alpha\,a)\,[\alpha\,(aa^{-1})] = [a\,(aa^{-1})]\,(aa^{-1})$. Therefore,

(2):   $$a^{-1} = [a^{-1}(a^{-1}a)]\,(a^{-1}a).$$

Then $[a^{-1}(a^{-1}a)]\,(ab) = [a^{-1}(a^{-1}a)]\ \{[(aa^{-1})a]\,b\} = [[\,\{\,[(a^{-1}(a^{-1}a))\,(aa^{-1})]\,a\}\,]]\,b = (a^{-1}a)\,b = (a^{-1}a)^2\,b$. So

(3):   $$[a^{-1}(a^{-1}a)]\,(ab) = (a^{-1}a)\,b = (a^{-1}a)^2\,b$$

By (3), $[a^{-1}(a^{-1}a)]\,(ab) = (a^{-1}a)\,b = (a^{-1}a)\,[(bb^{-1})^2\,b] = \{[(a^{-1}a)\,(\,b^{-1}b)]\,(bb^{-1})\}\,b = \{[(a^{-1}a)\,(b^{-1}b)]\}\,b$, so

(4): $\qquad (a^{-1} a)\, b = \{[(a^{-1} a)\, (b^{-1} b)]\}\, b.$

By (3) and (4), we get

(5): $[(b^{-1} b)\, (a^{-1} a)]\ \{[(a^{-1} (a^{-1} a)]\, (ab)\} = [(b^{-1} b)\, (a^{-1} a)]\, [(a^{-1} a)^2\, b] = [(b^{-1} b)\, (a^{-1} a)]\, b = [(a^{-1} a)\, (b^{-1} b)]\, b.$

Using (1), (4) and (5) we have that $\{(ab)\, [(\alpha\, b^{-1})\, (\alpha\, a^{-1})]\}\, (ab) = [(ab)b^{-1}]\ \{[(\,a^{-1} a)]\, (\,b^{-1} b)]\, b\} =$

$= a\ [[\ (bb^{-1})\ \{[(a^{-1} a)]\, (\,b^{-1} b)]\, b\}\ ]] =$

$=\ a\ \{[(a^{-1} a)\, (b^{-1} b)]\, b\} =\ a\ \{(a^{-1} a)\, b\} = a\ \{(a^{-1} a)^2\, b\} = \{[a(aa^{-1})]\, (aa^{-1})\}\, b.$

and so, by (2), $\{(ab)\, [(\alpha\, b^{-1})\, (\alpha\, a^{-1})]\}\, (ab) = ab.$

Note that $(a^{-1} a)\, [a^{-1} (a^{-1} a)\ ] = [\alpha\, (a^{-1} a)]\, (\alpha\, a^{-1}) = \alpha\, [(a^{-1} a)a^{-1}] = \alpha\, a^{-1} = a^{-1} (a^{-1} a\,)$ and so

(6): $\qquad\qquad\qquad (a^{-1} a)\, [a^{-1} (a^{-1} a)] = a^{-1} (a^{-1} a)$

Note that $(a^{-1} a)b = (a^{-1} a)\, [(b^{-1} b)^2\, b] = \{\ [(a^{-1} a)\, (b^{-1} b)]\, (b^{-1} b)\ \}\ b = [(a^{-1} a)\, (b^{-1} b)]b$ and so

(7): $\qquad\qquad\qquad (a^{-1} a)b = [(a^{-1} a)\, (b^{-1} b)]b\ = [(b^{-1} b)\, (a^{-1} a)]b$

We now wish to show that $\{[(\alpha\, b^{-1})\, (\alpha\, a^{-1})]\, (ab)\}\ \ [(\alpha\, b^{-1})\, (\alpha\, a^{-1})]\ =\ (\alpha\, b^{-1})\, (\alpha\, a^{-1}).$ Now

$\{\ [(\alpha\, b^{-1})\, (\alpha\, a^{-1})]\, (ab)\ \}\ \ \ [(\alpha\, b^{-1})\, (\alpha\, a^{-1})]\ =$

$=\qquad <\ \{[(b^{-1} (b^{-1} b)]\, [a^{-1} (a^{-1} a)]\}\ \ \{[(aa^{-1})a]\, b\}\ >\qquad\qquad \{\ [(b^{-1} (b^{-1} b))\, [a^{-1} (a^{-1} a)]\ \}\ \ =$

$=\quad <\ \{\ [[\ \{[(b^{-1} (b^{-1} b))\, [a^{-1} (a^{-1} a)]\}\ \ (aa^{-1})\ ]]\ a\ \}\ b\ >\qquad \{\ [(b^{-1} (b^{-1} b))\, [a^{-1} (a^{-1} a)]\ \}\ =$

$=\ <\ \{\{\ \ [(b^{-1} (b^{-1} b))\ \ [[\ \{\ [a^{-1} (a^{-1} a)]\, (aa^{-1})\ \}\ a\ ]]\ \}\}\ b\ >\ \ \{\ [(b^{-1} (b^{-1} b))\, [a^{-1} (a^{-1} a)]\ \}\ =$

which by (2) equals $\ [\ \{\ [(b^{-1} (b^{-1} b)]\, (a^{-1}\ a\ ))\ \}\ b\ ]\ \ \ \{\ [(b^{-1} (b^{-1} b))\, [a^{-1} (a^{-1} a)]\ \}\ =$

$=\ (\ b^{-1}\, \{\ [(b^{-1} b)\, (a^{-1} a)]\, b\ \}\ \ )\ \ \ \{\ [(b^{-1} (b^{-1} b))\, [a^{-1} (a^{-1} a)]\ \}$ , which by (7) equals

$=\ \{\ b^{-1}\, [(a^{-1} a)b]\ \}\ \ \ \{\ [(b^{-1} (b^{-1} b))\, [a^{-1} (a^{-1} a)]\ \}\ =$

$=\ <\ \ [\ \{\ b^{-1}\, [(a^{-1} a)b]\ \}\ b^{-1}]\ \ (b^{-1} b)\ >\ \ \ [a^{-1} (a^{-1} a)]\ =$

$=\ \{\ b^{-1}\ [\ \{[(a^{-1} a)b]\, b^{-1}\}\ (b^{-1} b)\ ]\ \}\ \ \ [a^{-1} (a^{-1} a)]\ =$

$= \ < b^{-1} \ \{(a^{-1}a)\,[(b\,b^{-1})(b^{-1}b)]\} \ > \qquad [a^{-1}(a^{-1}a)] \ =$

$= \ \{ \ b^{-1} \ [(b^{-1}b)^2 \ (a^{-1}a)] \ \} \qquad [a^{-1}(a^{-1}a)] \qquad =$

$= \ < \{[b^{-1}(b^{-1}b)]\,(b^{-1}b)\} \ (a^{-1}a) \ > \quad [a^{-1}(a^{-1}a)] \ =$

$= \ [b^{-1}(b^{-1}b)] \qquad \{ \ [(b^{-1}b)\,(a^{-1}a)] \ [a^{-1}(a^{-1}a)] \ \} \ =$

$= \ [b^{-1}(b^{-1}b)] \qquad < \ (b^{-1}b)^2 \ \{(a^{-1}a)^2 \ [a^{-1}(a^{-1}a)]\} \ > \ =$

$= \ < \ \{[b^{-1}(b^{-1}b)]\,(b^{-1}b)\} \ (b^{-1}b) \ > \quad \{(a^{-1}a)^2 \ [a^{-1}(a^{-1}a)]\} \ , \quad$ which by (2) and (6) equals

$= \ [b^{-1}(b^{-1}b)] \ [a^{-1}(a^{-1}a) \ = (\alpha\,b^{-1}) \ (\alpha\,a^{-1}).$

So we have proved that $(ab)^{-1} = (b^{-1}*b^{-1}b) \ (a^{-1}*a^{-1}a) = (\alpha\,b^{-1}) \ (\alpha\,a^{-1})$ for all $a,b \in S$ and so E($S$) is a semilattice **_and_** $(ab)^{-1} = (\alpha\,b^{-1}) \ (\alpha\,a^{-1})$ for all $a,b \in S$.

We now want to prove that $(ab)c = (\alpha\,a)(bc)$ for all $a,b,c \in S$. Recall that, from (2) and (5), $a = [a\,(aa^{-1})] \ (aa^{-1}) = (\alpha\,a) \ (aa^{-1})$ and $(aa^{-1}) \ (\alpha\,a) = \alpha\,a$. Now $(\alpha\,a) \ (bc) = (\alpha\,a) \ \{ \ [(\alpha\,b) \ (bb^{-1})] \ c \ \} =$
$= \{ \ [(\alpha\,a) \ (\alpha\,b)] \ (bb^{-1}) \ \} \ c = \ \{[(\alpha\,a) \ (\alpha\,b)] \ [\alpha\,(bb^{-1})] \ \} \ c = \{ \ \alpha\,[(ab) \ (\,bb^{-1})] \} \ c$. So we need only prove that $(ab)(bb^{-1}) = \alpha\,(ab)$, for then $(\alpha\,a)(bc) = [\alpha\,\alpha\,(ab)] \ c = (ab) \ c$. To prove this we now use the facts that E($S$) is a semilattice **_and_** $(ab)^{-1} = (\alpha\,b^{-1}) \ (\alpha\,a^{-1})$ for all $a,b \in S$. We have that

$[(ab)^{-1} \ (ab)] \ (bb^{-1}) = (bb^{-1}) \ [(ab)^{-1} \ (ab)] = \ (bb^{-1}) \ \{ \ [(\alpha\,b^{-1}) \ (\alpha\,a^{-1})] \ (ab) \ \} = \{[(bb^{-1}) \ (\alpha\,b^{-1})] \ (\alpha\,a^{-1})]\} \ (ab)$
$= [(\alpha\,b^{-1}) \ (\alpha\,a^{-1})] \ (ab) = (ab)^{-1} \ (ab)$. Therefore, $(ab)(\,bb^{-1}) = \ \{ \ [(ab) \ (ab)^{-1}] \ (ab) \ \} \ (bb^{-1}) =$
$= (ab) \ \{ \ [(ab)^{-1} \ (ab)] \ (bb^{-1}) \ \} = (ab) \ [(ab)^{-1} \ (ab)] = (ab)[(ab) \ (ab)^{-1}] = \alpha\,(ab)$. This completes the proof of the sufficiency of the conditions in Lemma 8, and therefore completes the proof of Lemma 8. ∎

**Lemma 9.** *$S$ is determined by a mapping $\alpha \in \text{AUT}_e^2(S,*)$ and $(S,*)$ is a semigroup that is a semilattice* E *of groups if and only if*

(9.1) *$S$ is a completely inverse, generalized right-Bol groupoid,*
(9.2) *$\alpha : a \mapsto a\,(aa^{-1})$ satisfies $\alpha \in \text{AUT}^2(S)$ and*
(9.3) E($S$) *is a semilattice **_or_** $(ab)^{-1} = (\alpha\,b^{-1}) \ (\alpha\,a^{-1})$ for all $a,b \in S$.*

Proof. ($\Longrightarrow$) This follows from Lemma 1.
($\Longleftarrow$) Assume that (9.1), (9.2) and (9.3) are valid. It follows from Lemma 8 that $(ab)c = (\alpha\,a)(bc)$ for all $a,b,c \in S$. Also, (9.2) implies that $\alpha\,e = e \ (ee^{-1}) = e$ and so $\alpha \in \text{AUT}_e^2(S)$. Then it follows from Lemma 7 that $S$ is determined by $\alpha \in \text{AUT}_e^2(S,*)$, with $(S,*)$ a semilattice of groups. This completes the proof of Lemma 9. ∎

The next Lemma establishes conditions under which (1.4) and (1.6) are equivalent.

**Lemma 10.** *If S is a completely inverse, generalized right-Bol groupoid, $\alpha \in \mathrm{AUT}_e^2(S)$ and either*

E(S) *is a semilattice* **or** $(ab)^{-1} = (\alpha\, b^{-1})\, (\alpha\, a^{-1})$ *for all a,b* $\in$ *S*

*then* $\alpha : a \mapsto a\, (aa^{-1})$ *for all* $a \in S$ *if and only if* $(ab)c = (\alpha\, a)(bc)$ *for all a,b,c* $\in$ *S.*

Proof. ($\Longrightarrow$) Lemma 9 implies that $S$ is determined by $\alpha \in \mathrm{AUT}_e^2(S,*)$ and the semigroup $(S,*)$ is a semilattice of groups. Then Lemma 1 implies that $(ab)c = (\alpha\, a)(bc)$ for all $a,b,c \in S$.

($\Longleftarrow$) Lemma 4 implies that E(S) is a semilattice. So Lemma 7 implies that $S$ is determined by $\alpha \in \mathrm{AUT}_e^2(S,*)$, with the semigroup $(S,*)$ a semilattice of groups. Then Lemma 1 implies that $\alpha : a \mapsto a\, (aa^{-1})$ for all $a \in S$. ∎

## 4. Semilattices of groupoids determined by groups

**Lemma 11.** *Let* E *be any semilattice, and to each e* $\in$ E *assign a groupoid S(e), where S(e) is determined by an involutive, idempotent-fixed mapping $\alpha_e$ on G(e), with each G(e) (e* $\in$ E*) a group and such that G(e) and G(f) are disjoint if e* $\neq$ *f in* E*. To each pair of elements f, e of* E *such that f* $\geq$ *e, assign a homomorphism $\delta_{f,e} : S(f) \rightarrow S(ef)$ such that if g* $\geq$ *f* $\geq$ *e then*

(A): $$\delta_{f,e}\, \delta_{g,f} = \delta_{g,e}.$$

*Assume that $\delta_{e,e}$ is the identity automorphism of S(e). Assume also that, for b* $\in$ *S(f) and all e* $\in$ E*,*

(B): $$\alpha_e(\delta_{f,e}\, b) = \delta_{f,e}(\alpha_f\, b).$$

*Let S be the union of all the S(e), (e* $\in$ E*), and define the product of any two elements a* $\in$ *S(e) and b* $\in$ *S(f) as*

(C): $$ab = (\delta_{e,ef}\, a)\, (\delta_{f,ef}\, b).\ Then$$

(11.1) *Each $\delta_{f,e}$ is a group homomorphism from G(f) to G(e);*

(11.2) *The union G of all the groups G(e), (e* $\in$ E*), is a semilattice* E *of groups if we define a product $*$ of any two elements a* $\in$ *G(e) and b* $\in$ *G(f) as*

(D): $$a*b = (\delta_{e,ef}\, a) \circ (\delta_{f,ef}\, b),\ where\ the\ product\ \circ\ is\ the\ product\ in\ the\ group\ G(ef);$$

(11.3) *S is determined by the mapping $\alpha \in \mathrm{AUT}_e^2(G,*)$, where $\alpha = \bigcup_{e \in E} \alpha_e$ ( that is, $\alpha\, a = \alpha_e a\ (a \in G(e)\,)$*

*and (G,$*$) is a semilattice* E *of groups.*

Proof. Let $f \geq e$ for some $e,f \in$ E. By definition, $G(f) = S(f)$, $G(ef) = S(ef) = G(e)$ and $\delta_{f,e} : S(f) \rightarrow S(ef)$. Hence, $\delta_{f,e} : G(f) \rightarrow G(ef) = G(e)$. Let $+$ denote the product in $G(f)$ and let $\circ$ denote the product in $G(ef) = G(e)$ and let $a,b \in G(f)$. Then $\delta_{f,e}(a+b) = \delta_{f,e}[(\alpha_f\, a)b] = [\delta_{f,e}(\alpha_f\, a)](\delta_{f,e}\, b) = \{\alpha_e[\delta_{f,e}(\alpha_f\, a)]\} \circ (\delta_{f,e}\, b)$. Using (B), $\delta_{f,e}(a+b) = \{\alpha_e[\,\delta_{f,e}(\alpha_f\,a)]\} \circ (\delta_{f,e}\,b) = \delta_{f,e}[\,\alpha_f(\alpha_f\,a)] \circ (\delta_{f,e}\,b) = (\delta_{f,e}\,a) \circ (\delta_{f,e}\,b)$ and so

$\delta_{f,e}: G(f) \to G(e)$ is a homomorphism. This proves (11.1).

Now (11.2) is a consequence of (11.1), (A) and [1, Theorem 4.11].

Let $a \in S(e)$ and $b \in S(f)$, for any $e, f \in E$. Let $\circ$ denote the product in $G(ef)$. Then, by the definition of the product in $S$ and using (B), $ab = (\delta_{e,ef}\,a)(\delta_{f,ef}\,b) = [\alpha_{ef}(\delta_{e,ef}\,a)] \circ (\delta_{f,ef}\,b) = [\delta_{e,ef}(\alpha_{ef}\,a)] \circ (\delta_{f,ef}\,b) = (\alpha_{ef}\,a)*b =$
$= (\alpha\,a)*b$. We need only prove that $\alpha \in \mathrm{AUT}_e^2(G,*)$. We have that $\alpha(a*b) = \alpha[(\delta_{e,ef}\,a) \circ (\delta_{f,ef}\,b)] =$
$\alpha_{ef}[(\delta_{e,ef}\,a) \circ (\delta_{f,ef}\,b)] = [\alpha_{ef}(\delta_{e,ef}\,a)] \circ [\alpha_{ef}(\delta_{f,ef}\,b)] = [\delta_{e,ef}(\alpha_{ef}\,a)] \circ [\delta_{f,ef}(\alpha_f\,b))] =$
$= [\delta_{e,ef}(\alpha\,a)] \circ [\delta_{f,ef}(\alpha\,b))] = (\alpha\,a)*(\alpha\,b)$ and so $\alpha$ is a homomorphism. Finally, $\alpha\,\alpha\,a = \alpha(\alpha_e a) = \alpha_e(\alpha_e a)$
$= a$ and so $\alpha$ is an automorphism. Clearly it is E-fixed. Hence, $\alpha \in \mathrm{AUT}_e^2(G)$. This proves (11.3) and completes the proof of Lemma 11.

Note that in Lemma 11 if $a = aa$ in $S$, where $S$ has the product defined in (C), then $a$ is an identity element of a group $G(e)$ for some $e \in E$. Conversely, every identity of a group $G(e)\,(e \in E)$ is an idempotent of the groupoid $S$. Hence, there is a natural bijection between E and $E(S)$ under which, in fact, $E \cong E(S)$. Therefore, we have Proved the first part of Corollary 12.

**Corollary 12.** *The groupoid $S$ in Lemma 11, with product defined as in* (C)*, is a semilattice of groupoids $S(e)$ $(e \in E(S))$, each of which is determined by an involutive, idempotent-fixed mapping $\alpha_e$ on a group $G(e)$. $S$ itself is determined by the involutive, idempotent-fixed mapping $\alpha = \bigcup_{e \in E} \alpha_e$ on the semigroup $(G,*)$ of Lemma 11, which is a semilattice of groups. $S$ satisfies the identity $(ab)c = (\alpha\,a)(bc)$.*

Proof. The second part of the Corollary follows from (11.3) and Result 2.

**Lemma 13.** *If $S$ is determined by a mapping $\alpha \in \mathrm{AUT}_e^2(S, \otimes)$ and $\{S, \otimes\}$ is a semigroup that is a semilattice* E *of groups then $S$ can be constructed as in Lemma 11 above.*

Proof: Let the groupoid $S$ be determined by a mapping $\alpha \in \mathrm{AUT}_e^2(S, \otimes)$, with $(S, \otimes)$ a semigroup and a semilattice E of groups. So (1.1) through (1.13) apply. Recall that we can consider $E = E(S)$.

First we define, for $e \in E$, $S(e) = \{a \in S : aa^{-1} = e\}$ where juxtaposition denotes product in the groupoid $S$ and where, as in the proof of (MT1), $a^{-1} = \alpha\,\mathrm{a}^{-1}$, where $\mathrm{a}^{-1}$ is the (unique) inverse of $a$ in $(S, \otimes)$. Define $G(e)$ to be the maximal subgroup of $(S, \otimes)$, with identity element $e$ [cf. 1, Theorem 1.11]. Thus, $G(e) = \{a \in S : a \otimes \mathrm{a}^{-1} = e\}$. By (1.13), $S(e) = G(e)$, for $e \in E$.

We show that $S(e)$ is determined by $(G(e), \otimes)$ and $\alpha_e = \alpha/_{G(e)}$, the restriction of $\alpha$ to $G(e)$. Since, for $a, b \in S(e) = G(e)$, $ab = (\alpha\,a) \otimes b = (\alpha_e a) \otimes b$. Clearly, $\alpha_e \in \mathrm{AUT}_e^2\,G(e)$. Therefore

(1): Each $S(e)$ is determined by $(G(e), \otimes)$ and $\alpha_e = \alpha/_{G(e)} \in \mathrm{AUT}_e^2\,G(e)$.

Note that since for $a \in S(e)$, $(\alpha\,a)(\alpha\,a)^{-1} = (\alpha\,a)(\alpha\,a^{-1}) = \alpha(aa^{-1}) = aa^{-1} = e$, $\alpha\,a \in S(e)$ and so $\alpha\, S(e) \subseteq S(e)$. Also, $a = \alpha\,\,\alpha\,\,a \in \alpha\, S(e)$, so $\alpha\, S(e) \subseteq S(e) \subseteq \alpha\, S(e)$ and hence we have proved that

(2):  $\qquad\qquad\qquad\qquad S(e) = \alpha\, S(e)$. Also $a \in S(e)$ if and only if $\alpha\, a \in S(e)$.

Note also that for $a \in S(e)$, $ea = (aa^{-1})\, a = a$ and , by (1.10), $ae = e\, (\alpha\, a) = \alpha\, (ea) = \alpha\, (a)$. We have proved

(3):  $\qquad\qquad\qquad\qquad$ if $a \in S(e)$ then $ea = a$ and $ae = \alpha\, a$.

Now let $a \in S(e)$ and $b \in S(f)$, for some $e, f \in$ E. Then, using (3), (2.6) and (2.10), we have that $(af)\,(be) = [f\,(\alpha\, a)]\,[e\,(\alpha\, b)] = (\alpha\, f)\,\{(\alpha\, a)\,[e\,(\alpha\, b)]\} = f\,\{(\alpha\, a))[e\,(\alpha\, b)]\} = f\,[(ae)\,(\alpha\, b)] = f\,[(\alpha\, a)\,(\alpha\, b)] = $ $= f\,[\,\alpha\,(ab)] = (ab)\,f = (\alpha\, a)\,(bf) = (\alpha\, a)\,(\alpha\, b) = \alpha\,(ab)$. Therefore,

$ab = \alpha\,\,\alpha\,(ab) = \alpha\,[(af)\,(be)] = [\,\alpha\,(af)]\,[\,\alpha\,(be)] = [(\alpha\, a)\,(\alpha\, f)][(\alpha\, b)\,(\alpha\, e)] = [(\alpha\, a)\,f]\,[(\alpha\, b)\, e] = (fa)\,(eb)$. We have proved that

(4):  $\qquad\qquad\qquad\qquad$ if $a \in S(e)$ and $b \in S(f)$, for some $e, f \in$ E then
$\qquad\qquad\qquad\qquad\qquad ab = (fa)\,(eb)$.

To each pair $\{f, e\} \subseteq$ E such that $f \geq e$, define $\delta_{f,\, e}: S(f) \rightarrow S(ef)$ as $\delta_{f,\, e}\, b = eb\ (b \in S(f))$. Then $\delta_{f,\, e}$ is well defined, as a consequence of (1.12), since $(eb)(eb)^{-1} = (ee^{-1})(bb^{-1}) = ef \in S(ef)$ .

Now if $a \in S(f)$ and $b \in S(f)$ then, by (1.11), for any $e \in$ E, $(ea)(eb) = e(ab)$ and so, $(\delta_{f,\, e}\, a)(\delta_{f,\, e}\, b) = (ea)(eb) = e(ab) = \delta_{f,\, e}\,(ab)$ . We have proved that

(5):  $\qquad\qquad\qquad\qquad$ Each $\delta_{f,\, e}: S(f) \rightarrow S(ef)$ is a homomorphism.

If $g \geq f \geq e$ $(e, f, g \in$ E) and $a \in S(g)$ then $\delta_{f,\, e}\,\delta_{g,\, f}\, a = \delta_{f,\, e}\,(fa) = e\,(fa) = [(\alpha\, e)\, f]\, a = (ef)\, a = ea = \delta_{g,\, e}\, a$ and so

(6):  $\qquad\qquad\qquad\qquad$ if $g \geq f \geq e$ $(e, f, g \in$ E) then $\delta_{f,\, e}\,\delta_{g,\, f} = \delta_{g,\, e}$ .

Also, from (3), if $a \in S(e)$ and $b \in S(f)$, for some $e, f \in$ E then $ab = (fa)(eb)$. But $(ef)a = (fe)a = (\alpha\, f)(ea) = fa$ and $(ef)b = (\alpha\, e)(fb) = eb$ and so, $ab = (\delta_{e,\, ef}\, a)(\delta_{f,\, ef}\, b)$; that is,

(7):  $\qquad\qquad\qquad\qquad$ if $a \in S(e)$ and $b \in S(f)$, for some $e, f \in$ E then
$\qquad\qquad\qquad\qquad\qquad ab = (\delta_{e,\, ef}\, a)(\delta_{f,\, ef}\, b)$.

In addition, if $b \in S(f)$ and $f \geq e$ for some $e, f \in$ E then $\alpha_e\,(\delta_{f,\, e}\, b) = \alpha_e\,(eb) = \alpha\,(eb) = (\alpha\, e)(\alpha\, b) = e\,(\alpha\, b) = $ $= \delta_{f,\, e}\,(\alpha\, b) = \delta_{f,\, e}\,(\alpha_f\, b)$, which proves that

(8):  $\qquad\qquad\qquad\qquad$ If $b \in S(f)$ and $f > e$ for some $e, f \in$ E then
$\qquad\qquad\qquad\qquad\qquad \alpha_e\,(\delta_{f,\, e}\, b) = \delta_{f,\, e}\,(\alpha_f\, b)$

Now, consider the groupoids $S(e)$ $(e \in$ E), each of which - by (1) - is determined by $\alpha_e = \alpha\,/_{G(e)} \in \mathrm{AUT}_e^2\ G(e)$. As a consequence of (5), (6), (7) and (8), we also have a collection of homomorphisms $\delta_{f,\, e}: S(f) \rightarrow S(ef)$ $(f \geq e$ in E) that satisfy (A), (B) and (C) of the first part of Theorem 11. Therefore, as in the proof of Theorem 11, (11.1), (11.2), (D) and (11.3) are valid. Hence, the union $S$ of the $S(e)$ ($e \in$ E), with product

$$ab = (\delta_{e, ef}\, a)(\delta_{f, ef}\, b) = [\,\alpha\,(\delta_{e, ef}\, a)] \circ (\delta_{f, ef}\, b) = [\,\delta_{e, ef}\,(\alpha_e\, a))] \circ (\delta_{f, ef}\, b) \ \ \text{for } a \in S(e) \text{ and } b \in S(f),$$

is determined by $\alpha \in \mathrm{AUT}_e^2 (G, *)$, with the semigroup $(G, *)$ a semilattice E of the groups $G(e)$ ($e \in$ E), and where $a * b = (\delta_{e, ef}\, a) \circ (\delta_{f, ef}\, b)$ and $\alpha = \bigcup_{e \in E} \alpha_e$. That is, $ab = (\alpha_e\, a) * b = (\alpha\, a) * b$.

Now, by the hypothesis of Lemma 13, $a \otimes b = (\alpha\, a)b = a * b$. This assures that the groupoid $S$ in Lemma 12 can be constructed as required. ∎

**Lemma 14.** *If $S$ is a semilattice* E$(S)$ *of groupoids $S(e)$ ($e \in$ E), where each $S(e)$ is determined by a mapping $\alpha_e \in \mathrm{AUT}_e^2\, G(e)$ and each $G(e)$ ($e \in$ E) is a group with identity $e$, then $S$ is a completely inverse groupoid.*

Proof. Let $e \in$ E$(S)$ and $a \in S(e)$. Then, as in the proof of Theorem 2, $x = \alpha_e\, a^{-1}$ is the unique inverse of $a$ in the completely inverse groupoid $S(e)$, where $a^{-1}$ is the unique inverse of $a$ in the group $G(e)$. If $y \in S(f)$ is an inverse of $a$ in $S$ then, since $S$ is a semilattice E$(S)$ of groupoids $S(e)$ ($e \in$ E), $\{a, y\} \subseteq S(ef) = S(e)$ and so $y = x$. Then, since by Lemma 1 each $S(e)$ is completely inverse, so is $S$. ∎

## 5. Proof of the Main Theorem.

Proof. As a consequence of Lemma 1, (M1) implies (M2), (M3) and (M4). By Lemma 5, (M2) implies (M1). By Lemma 7, (M3) implies (M.1). By Lemma 9, (M4) implies (M1). Hence, conditions (M.1), (M2), (M3) and (M4) are equivalent.

By Lemmas 11 and 13, (M1) and (M5) are equivalent. By Lemma 13 and the proof of Lemma 11 and Corollary 12, (M1) implies (M6).

We now prove that (M6) implies (M1). Using Lemma 14, $S$ in (M6) is completely inverse. Since (M1) and (M2) are equivalent statements, we need only prove that $\alpha \in \mathrm{AUT}^2 (S)$, as the hypotheses of (M6) state that E$(S)$ is a semilattice. Firstly, we prove that $\alpha\, e = e$, for all $e \in$ E$(S)$. Since $e = (ee)e = (\alpha\, e)(ee) = (\alpha\, e)e$ and since $\alpha\, e = \alpha_e\, e \in S(e)$, $\alpha\, e = \alpha_e\, e = \alpha_e\,[(\alpha\, e)e] = \alpha_e\,[(\alpha_e\, e)e] = [\alpha_e\,(\alpha_e\, e)](\alpha_e\, e) = e(\alpha\, e) = [(\alpha\, e)e]\,(\alpha\, e)$ and $[e(\alpha\, e)]e = (\alpha\, e)e = e$. So

(1): $$\alpha\, e = e^{-1} = e.$$

Then for any $a \in S(e)$, $a = (aa^{-1})a = (\alpha\, a)(a^{-1}a) = (\alpha\, a)(e) = (\alpha\, a)(ee) = (ae)e$ and so

(2): $$a = (\alpha\, a)e = (ae)e \ \ , \text{ for any } a \in S(e)$$

Then from (2), $\alpha\, a = \alpha_e\,[(ae)e] = [\alpha_e\,(ae)](\alpha_e\, e) = [(\alpha_e\, a)(\alpha_e\, e)]e = [(\alpha_e\, a)e]e = [(\alpha\, a)e]e = ae$ and so

(3): $$\alpha\, a = ae.$$

Now let $a \in S(e)$ and $b \in S(f)$ for any $e, f \in \mathrm{E}(S)$. Then, using (2) and (3), $(\alpha\, a)(\alpha\, b) = (\alpha\, a)(bf) = (ab)f =$
$= [\alpha_{ef}\, \alpha_{ef}\, (ab)]f = \{[(ab)(ef)](ef)\}f = \{[\alpha\,(ab)](ef)\}f = (ab)\, [(ef)f] = (ab)(ef) = \alpha\,(ab)$, proving that $\alpha \in \mathrm{AUT}^2(S)$.

Hence, (M6) implies (M1) and the proof of the Main Theorem is complete. ■

The author is grateful for the referee's close scrutiny of the original version of this paper and for suggestions that greatly improved its presentation. He also thanks the editors for their abundant patience and assistance.